\documentclass{jfm}
\usepackage{xspace,colortbl}
\usepackage[square,comma]{natbib}
\usepackage{amssymb}
\usepackage{epsfig}
\usepackage{graphicx}
\usepackage{amsmath}
\usepackage{bm}
\input{tcilatex}

\begin{document}

\title[Archimedes' famous-theorem]
{Archimedes' famous-theorem}

\author[H. Gouin ]{Henri Gouin\footnote{dedicato ai miei amici italiani}} \affiliation {Aix-Marseille Universit\'e, CNRS, Centrale Marseille, M2P2 UMR 7340  \\     13451  Marseille  France

  E-mails:\,\ henri.gouin@univ-amu.fr;\,\ henri.gouin@yahoo.fr
 \\

}

\maketitle

\begin{abstract}

In his treatise addressed to Dositheus of Pelusium (\cite{Calinger}), Archimedes of Syracuse obtained the result of which he was the most proud: \emph{a sphere has two-thirds the volume  of its circumscribing cylinder} (see Fig. 1).    At his request a sculpted sphere and cylinder were placed on his tomb near Syracuse (\cite{Rorres}).\\
\noindent Usually, it is admitted that to find this formula, Archimedes used a half polygon inscribed in a semicircle; then he performed rotations of these two figures to obtain a set of trunks in a sphere. This set of trunks  allowed him to determine the volume.\\
\noindent In our opinion, Archimedes was so clever that he found the proof with shorter demonstration. Archimedes did not need to know $\pi$ to prove the result and the Pythagorean theorem was probably the key to the proof.

\end{abstract}


\textbf{----------------------------------------------------------------------------------------------- }

{\qquad\qquad\qquad\qquad\qquad\quad\quad\qquad\qquad\emph{\textbf{Proof}}}

\noindent \emph{See Fig. 2a:}

The area generated by the rotation of segment $AC$ around $\overrightarrow{Oz}$ at level $\ell$  is proportional to the square of radius $r$  (the proportionality coefficient is $\pi$).

The area generated by the rotation of segment $AB$ around $\overrightarrow{Oz}$ at level $\ell$  is proportional to the radius square of a radius which is, on the basis on the Pythagorean theorem, $\sqrt{r^ 2 - \ell^2}$.

The area generated by the rotation of segment $BC$ around $\overrightarrow{Oz}$ at level $\ell$  is proportional to $r^ 2 - (r^2 -\ell^2) =\ell^2$.\\

\noindent \emph{See Fig. 2b:}

Volume $V_1$ of the domain included between the cylinder and  the sphere is therefore twice the volume of a cone whose the base is that of the cylinder and the height is the radius of the sphere.
\\
(Archimedes  knew the volume of the cone which is the extension of the volume of a  pyramid with $n$ faces, $(3 \leq n \in N)$).\\

\noindent \emph{Consequence:}

Two times the volume of the cone is $V_1 =  (2/3)\, r\, S$ where S is the surface of the cone basis.

The volume of the cylinder is $ V_2 =  2r  \, S $.

The volume of the sphere is $V_3 = V_2 - V_1 $ \\

\noindent \emph{Conclusion:}

$$
V_3 = \frac{2}{3}\, V_2
$$
$\qquad\qquad\qquad\qquad\qquad\qquad\qquad\qquad\qquad\qquad\qquad\qquad\qquad\qquad\qquad\qquad\qquad\qquad\qquad\square$

\noindent \emph{\textbf{ {Remark:} }}
If we know - as Archimedes did (\cite{McKeeman}) -  the signification of $\pi$,
$$
 V_2  = 2 r\, ( \pi r^2 )\quad  {\rm and\ consequently} \quad V_3 =\frac{4}{3}\,\pi r^3\,.
$$
\begin{figure}
\begin{center}
\includegraphics[width=5cm]{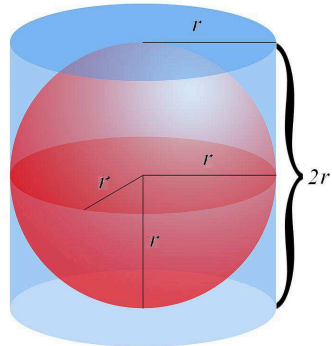}
\end{center}
\caption{\emph{A sphere has 2/3 the volume  of its circumscribing cylinder.\newline  The radius of the sphere is $r$.}} \label{fig1}
\end{figure}

\begin{figure}
\begin{center}
\includegraphics[width=11cm]{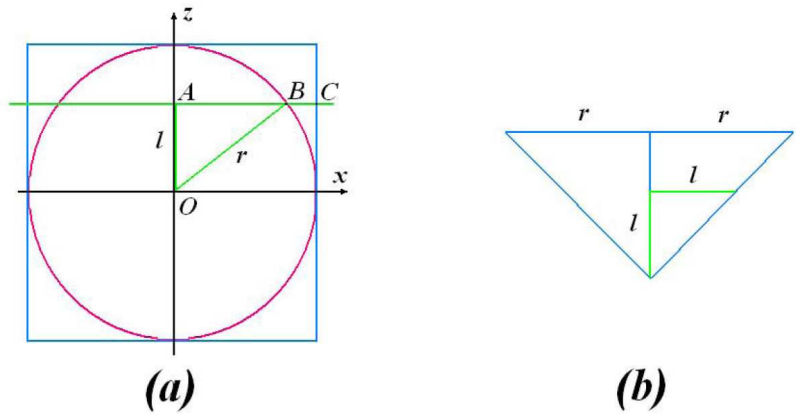}
\end{center}
\caption{ {Sphere, cylinder and cone.}} \label{fig1}
\end{figure}

\end{document}